\documentclass[11pt]{article}
\usepackage{amsmath}
\usepackage{amssymb}
\usepackage{amsfonts}
\usepackage{hyperref}
\usepackage[dvips]{graphicx}

\newtheorem{theorem}{Theorem}[section]

\newtheorem{lemma}[theorem]{Lemma}

\newtheorem{proposition}[theorem]{Proposition}

\let\originalleft\left
\let\originalright\right
\renewcommand{\left}{\mathopen{}\mathclose\bgroup\originalleft}
\renewcommand{\right}{\aftergroup\egroup\originalright}
\begin{document}

\title{Flows of co-closed $G_{2}$-structures}
\author{Sergey Grigorian \thanks{%
School of Mathematical and Statistical Sciences, University of Texas Rio
Grande Valley, Edinburg, TX 78539, USA} }
\maketitle

\begin{abstract}
We survey recent progress in the study of $G_{2}$-structure Laplacian
coflows, that is, heat flows of co-closed $G_{2}$-structures. We introduce
the properties of the original Laplacian coflow of $G_{2}$-structures as
well as the modified coflow, reviewing short-time existence and uniqueness
results for the modified coflow and well as recent Shi-type estimates that
apply to a more general class of $G_{2}$-structure flows.
\end{abstract}

\section{Introduction}

\setcounter{equation}{0}One of the most successful techniques in geometric
analysis has been the application of geometric flows to various problems in
geometry and topology, most notably the use of the Ricci flow \cite%
{Hamilton3folds, MorganTianRicci} to solve the Poincar\'{e} Conjecture \cite%
{PerelmanRicci}. The Ricci flow is a non-linear weakly parabolic partial
differential equation for the Riemannian metric $g$ 
\begin{equation}
\frac{\partial g}{\partial t}=-2\mathrm{Ric}_{g}  \label{Ricflow0}
\end{equation}%
so that the evolution of the metric is given by the Ricci curvature defined
by the metric. This can further be interpreted as a heat equation for the
metric. In $G_{2}$-geometry, there have been a number of proposals for
geometric flows of $G_{2}$-structures. The general idea is that given an
initial $G_{2}$-structure with weaker assumptions than vanishing torsion,
the flow should eventually seek out a torsion-free $G_{2}$-structure, if one
exists on the given manifold. A $G_{2}$-structure is defined by a positive $%
3 $-form $\varphi $, which in turn defines the metric $g$, and the
corresponding Hodge dual $4$-form $\ast \varphi =:\psi $. Therefore, a
natural equation to consider is the analog of the heat equation for the $3$%
-form $\varphi $%
\begin{equation}
\frac{\partial \varphi }{\partial t}=\Delta _{\varphi }\varphi .
\label{lapflowphi}
\end{equation}%
This Laplacian flow of the $3$-form $\varphi $ is now nonlinear in $\varphi $%
, because the metric and hence the Laplacian depend on $\varphi $ itself. A
particular case of this flow has been first studied by Robert Bryant \cite%
{bryant-2003}, where he restricted it to closed $G_{2}$-structures, that is
ones where $d\varphi =0$. For a closed $G_{2}$-structure, $\Delta \varphi
=dd^{\ast }\varphi $, so in this case, the $3$-form $\varphi $ stays closed
under the flow (\ref{lapflowphi}), and in fact remains within the same
cohomology class since $\Delta \varphi $ is exact. Short-time existence and
uniqueness of solutions to (\ref{lapflowphi}) was proved in \cite{BryantXu}.
Moreover, on a compact manifold $M$, this flow can be interpreted as the
gradient flow of the Hitchin functional $V$ given by%
\begin{equation}
V\left( \varphi \right) =\frac{1}{7}\int_{M}\varphi \wedge \ast _{\varphi
}\varphi .
\end{equation}%
The functional $V$ is then the volume of the manifold $M$. It was shown by
Nigel Hitchin in \cite{Hitchin:2000jd-arxiv} that if $\varphi $ is closed,
then the critical points of the functional $V$ within the cohomology class $%
\left[ \varphi \right] $ correspond precisely to torsion-free $G_{2}$%
-structures, and in particular, these critical points are maxima in the
directions transverse to diffeomorphisms. Under the flow (\ref{lapflowphi}), 
$V$ increases monotonically, so if the growth of $V$ is bounded, then $%
\varphi \left( t\right) $ would be expected to approach a torsion-free $%
G_{2} $-structure as $t\longrightarrow \infty $. The stability and
analyticity of this flow has recently been proved by Lotay and Wei \cite%
{LotayWei1,LotayWei2,LotayWei1a}

Alternatively, a $G_{2}$-structure and the corresponding metric may also be
defined by the $4$-form $\psi $ (up to a choice of orientation). Therefore,
instead of deforming $\varphi $, we may deform $\psi $. Using this idea,
Karigiannis, McKay, and Tsui, introduced the \emph{Laplacian coflow} for the 
$4$-form $\psi $ in \cite{KarigiannisMcKayTsui}. Instead of considering the
heat flow equation for $\varphi $, they instead considered the flow:%
\begin{equation}
\frac{\partial \psi }{\partial t}=\Delta _{\psi }\psi .  \label{lapcoflow}
\end{equation}%
If restricted to \emph{co-closed} $G_{2}$-structures (that is, ones with $%
d\psi =0$ and equivalently, those with a symmetric torsion tensor $T$) this
flow preserves the co-closed condition and in fact preserves the cohomology
class of $\psi $. In \cite{GrigorianCoflow}, it was shown that this flow has
similar characteristics to the original Laplacian flow for closed $G_{2}$%
-structures. In fact, (\ref{lapcoflow}) can also be regarded as a gradient
flow of the Hitchin functional (but now reformulated via $4$-forms).
However, a major difference compared with the Laplacian flow of closed $%
G_{2} $-structures (\ref{lapflowphi}) is that (\ref{lapcoflow}) is not even
a weakly parabolic equation. In fact, the symbol of the linearized equation
is indefinite. In order to have any hope of proving the existence of
solutions, a \emph{modified Laplacian coflow }of co-closed $G_{2}$%
-structures was introduced in \cite{GrigorianCoflow}:%
\begin{equation}
\frac{d\psi }{dt}=\Delta _{\psi }\psi +2d\left( \left( A-\mathop{\rm Tr}%
\nolimits T\right) \varphi \right)  \label{mod-coflow}
\end{equation}%
where $\mathop{\rm Tr}\nolimits T$ is the trace of the full torsion tensor $%
T $ of the $G_{2}$-structure defined by $\psi $, and $A$ is a positive
constant. This flow is now weakly parabolic in the direction of closed forms
and hence it is possible to relate it to a strictly parabolic flow using an
application of DeTurck's trick. Recently, the methods of Lotay and Wei for
Shi-type estimates for the flow (\ref{lapflowphi}) have been extended by Gao
Chen \cite{GaoChen1} to cover a more general class of $G_{2}$-structure
flows that includes (\ref{mod-coflow}) as well. We will first survey the
properties of $G_{2}$-structures and the Laplacian $\Delta _{\varphi
}\varphi $ in sections \ref{SecLapPsi} and \ref{SecDeformPhi}. Then, in
section \ref{secLapcoflow} we will focus on Laplacian coflows.

Despite the apparent similarity between closed and co-closed $G_{2}$%
-structures, there are also important differences. As shown in \cite%
{CrowleyNordstrom1}, co-closed $G_{2}$-structures always satisfy the $h$%
-principle (on both open and closed manifolds) and hence always exist
whenever a manifold admits $G_{2}$-structures. This is in contrast to closed 
$G_{2}$-structures for which the $h$-principle only holds on open manifolds.
Therefore, co-closed $G_{2}$-structures are in some sense more generic than
closed ones. This is both good and bad - it's good because they always
exist, but bad because one cannot expect their flows to always behave
nicely. This is also in part shown by the non-parabolicity of the original
coflow (\ref{lapcoflow}).

In this survey we will focus on analytic properties of flows on general $7$%
-manifolds, however another approach to understand the specific behavior of
geometric flows and obtain explicit solutions has been to consider manifolds
with some symmetry, in which case the number of degrees of freedom in the
PDE will be reduced. Both the original Laplacian coflow (\ref{lapcoflow})
and the modified Laplacian coflow (\ref{mod-coflow}) have been studied on a
variety of such manifolds with symmetry. Note that while in these situations
mostly the original coflow (\ref{lapcoflow}) with the negative sign has been
studied, results for the coflow with the positive sign (\ref{lapcoflow})
would be similar because equations reduce to ODEs. In \cite%
{KarigiannisMcKayTsui} and \cite{GrigorianSU3flow}, the coflow and the
modified coflow, respectively, have been studied on warped product manifolds
of the form $N^{6}\times L$ where $N^{6}$ is a $6$-dimensional manifold with 
$SU\left( 3\right) $-structure such as a Calabi-Yau or nearly K\"{a}hler
manifold and $L$ is either $\mathbb{R}$ or $S^{1}$. In particular, soliton
solutions in both cases have been obtained. In \cite{BagagliniFernandezFino1}%
, Bagaglini, Fernandez, and Fino, also studied both the coflows on the $7$%
-dimensional Heisenberg group. In particular, they have shown that the
long-term existence properties of the flow (\ref{mod-coflow}) depend on the
constant $A$. Similarly, in \cite{BagagliniFino1}, Bagaglini and Fino
studied the Laplacian coflow on $7$-dimensional almost-abelian Lie groups
and showed long-term existence properties and constructed soliton solutions.
In \cite{ManeroOtalVillacampa}, Manero, Otal, and Villacampa studied both
the Laplacian flow (\ref{cofloworig}) and the coflow (\ref{lapcoflow}) on
solvmanifolds, but instead of restricting to closed or co-closed $G_{2}$%
-structures, they instead restricted to \emph{locally conformally parallel }$%
G_{2}$-structures, which are the ones where only the $7$-dimensional $\tau
_{1}$ component of the torsion may be nonvanishing.

\subsubsection*{Acknowledgements}

The author is supported by the National Science Foundation grant DMS-1811754.

\section{Laplacian of a $G_{2}$-structure}

\label{SecLapPsi}\setcounter{equation}{0}Suppose $M$ is a smooth $7$%
-dimensional manifold with a $G_{2}$-structure $\varphi .$ Then we know $%
\varphi $ uniquely defines a compatible Riemannian metric $g_{\varphi }$,
the volume form $\mathop{\rm vol}\nolimits_{\varphi }$, Hodge star $\ast
_{\varphi }$, and the dual $4$-form $\psi =\ast _{\varphi }\varphi $. There
is arbitrary choice of orientation, which affects the relative sign of $\psi 
$. We use the same convention as \cite{Bryant-1987} and \cite%
{GrigorianG2TorsionWarped,GrigorianCoflow,GrigorianG2Torsion1,GrigorianSU3flow,GrigorianYau1}%
, which is opposite from the convention used in \cite%
{karigiannis-2005-57,karigiannis-2007}. For further properties of $\varphi $
and $\psi ,$ as well as different identities that they satisfy, we refer the
reader to the above references. We will also use the following notation. The
symbol $\lrcorner $ will denote contraction of a vector with the
differential form: 
\begin{equation}
\left( u\lrcorner \varphi \right) _{mn}=u^{a}\varphi _{amn}.
\end{equation}%
Note that we will also use this symbol for contractions of differential
forms using the metric, for example $\left( T\lrcorner \varphi \right)
_{a}=T^{mn}\varphi _{mna}$. Given a symmetric $2$-tensor $h$ on $M$, we
define the map $\mathrm{i}_{\varphi }:\Gamma \left( \mathop{\rm Sym}\nolimits%
\left( T^{\ast }M\right) \right) \longrightarrow \Lambda _{1}^{3}\oplus
\Lambda _{27}^{3}$ as 
\begin{equation*}
\mathrm{i}_{\varphi }\left( h\right) _{abc}=h_{[a}^{d}\varphi _{bc]d}^{{}}.
\end{equation*}%
We will define the operators $\pi _{1}$, $\pi _{7}$, $\pi _{14}$ and $\pi
_{27}$ to be the projections of differential forms onto the corresponding
representations. Sometimes we will also use $\pi _{1\oplus 27}$ to denote
the projection of $3$-forms or $4$-forms into $\Lambda _{1}^{3}\oplus
\Lambda _{27}^{3}$ or $\Lambda _{1}^{4}\oplus \Lambda _{27}^{4}$
respectively. For convenience, when writing out projections of forms, we
will sometimes just give the vector that defines the $7$-dimensional
component, the function that defines the $1$-dimensional component or the
symmetric $2$-tensor that defines the $1\oplus 27$ component whenever there
is no ambiguity. For instance,%
\begin{equation}
\begin{tabular}{lll}
$\pi _{1}\left( f\varphi \right) =f$ & $\pi _{1}\left( f\psi \right) =f$ & 
\\ 
$\pi _{7}\left( X\lrcorner \varphi \right) ^{a}=X^{a}$ & $\pi _{7}\left(
X\lrcorner \psi \right) ^{a}=X^{a}$ & $\pi _{7}\left( X\wedge \varphi
\right) ^{a}=X^{a}$ \\ 
$\pi _{1\oplus 27}\left( \mathrm{i}_{\varphi }\left( h\right) \right)
_{ab}=h_{ab}$ & $\pi _{1\oplus 27}\left( \ast \mathrm{i}_{\varphi }\left(
h\right) \right) _{ab}=h_{ab}$ & 
\end{tabular}
\label{projections}
\end{equation}%
The above-mentioned references give more information regarding the
properties of decomposition of differential forms with respect to $G_{2}$
representations.

The \emph{intrinsic torsion }of a $G_{2}$-structure is defined by $\nabla
\varphi $, where $\nabla $ is the Levi-Civita connection for the metric $g$
that is defined by $\varphi $. Following \cite{karigiannis-2007}, we have 
\begin{subequations}
\begin{eqnarray}
\nabla _{a}\varphi _{bcd} &=&T_{a}^{\ e}\psi _{ebcd}^{{}}  \label{codiffphi}
\\
\nabla _{a}\psi _{bcde} &=&-4T_{a[b}\varphi _{cde]}  \label{psitorsion}
\end{eqnarray}%
\end{subequations}%
where $T_{ab}$ is the \emph{full torsion tensor}. In general we can split $%
T_{ab}$ according to representations of $G_{2}$ into \emph{torsion components%
}: 
\begin{equation}
T=\frac{1}{4}\tau _{0}g-\tau _{1}\lrcorner \varphi +\frac{1}{2}\tau _{2}-%
\frac{1}{3}\tau _{3}  \label{Tdecomp}
\end{equation}%
where $\tau _{0}$ is a function, and gives the $\mathbf{1}$ component of $T$%
. We also have $\tau _{1}$, which is a $1$-form and hence gives the $\mathbf{%
7}$ component, and, $\tau _{2}\in \Lambda _{14}^{2}$ gives the $\mathbf{14}$
component and $\tau _{3}$ is traceless symmetric, giving the $\mathbf{27}$
component. As shown by Karigiannis in \cite{karigiannis-2007}, the torsion
components $\tau _{i}$ relate directly to the expression for $d\varphi $ and 
$d\psi $. In fact, in our notation, 
\begin{subequations}%
\label{dptors} 
\begin{eqnarray}
d\varphi &=&\tau _{0}\psi +3\tau _{1}\wedge \varphi +\ast \mathrm{i}%
_{\varphi }\left( \tau _{3}\right)  \label{dphi} \\
d\psi &=&4\tau _{1}\wedge \psi +\ast \tau _{2}.  \label{dpsi}
\end{eqnarray}%
\end{subequations}%
Note that in \cite{GrigorianCoflow,GrigorianG2Torsion1, GrigorianSU3flow,
GrigorianYau1} a different convention is used: $\tau _{1}$ in that
convention corresponds to $\frac{1}{4}\tau _{0}$ here, $\tau _{7}$
corresponds to $-\tau _{1}$ here,\ $\mathrm{i}_{\varphi }\left( \tau
_{27}\right) $ corresponds to $-\frac{1}{3}\tau _{3},$ and $\tau _{14}$
corresponds to $\frac{1}{2}\tau _{2}$. The notation used here is widely used
elsewhere in the literature.

An important special case is when the $G_{2}$-structure is said to be
torsion-free, that is, $T=0$. This is equivalent to $\nabla \varphi =0$ and
also equivalent, by Fern\'{a}ndez and Gray \cite{FernandezGray}, to $%
d\varphi =d\psi =0$. Moreover, a $G_{2}$-structure is torsion-free if and
only if the holonomy of the corresponding metric is contained in $G_{2}$ 
\cite{Joycebook}. On a compact manifold, the holonomy group is then
precisely equal to $G_{2}$ if and only if the fundamental group $\pi _{1}$
is finite. If $d\varphi =0$, then we say $\varphi $ defines a \emph{closed }$%
G_{2}$-structure. In that case, $\tau _{0}=\tau _{1}=\tau _{3}=0$ and only $%
\tau _{2}$ is in general non-zero. In this case, $T=-\frac{1}{2}\tau _{2}$
and is hence skew-symmetric. If instead, $d\psi =0$, then we say that we
have a \emph{co-closed }$G_{2}$-structure. In this case, $\tau _{1}$ and $%
\tau _{2}$ vanish in (\ref{dpsi}) and we are left with $\tau _{0}$ and $\tau
_{3}$ components. In particular, the torsion tensor $T_{ab}$ is now
symmetric.\ 

We will be using the following notation, as in \cite{GrigorianCoflow}. Given
a tensor $\omega $, the rough Laplacian is defined by 
\begin{equation}
\Delta \omega =g^{ab}\nabla _{a}\nabla _{b}\omega =-\nabla ^{\ast }\nabla
\omega .  \label{roughlap}
\end{equation}%
whereas the Hodge Laplacian defined by $\varphi $ or $\psi $ will be denoted
by $\Delta _{\varphi }$ or $\Delta _{\psi }$, respectively. For a vector
field $X$, define the \emph{divergence} of $X$ as 
\begin{equation}
\mathop{\rm div}\nolimits X=\nabla _{a}X^{a}.  \label{divX}
\end{equation}%
This operator can be extended to a $2$-tensor $\beta $: 
\begin{equation}
\left( \mathop{\rm div}\nolimits\beta \right) _{b}=\nabla ^{a}\beta _{ab}.
\label{divh}
\end{equation}%
Also, for a vector $X$, we can use the $G_{2}$-structure $3$-form $\varphi $
to define a \textquotedblleft curl\textquotedblright\ operator, similar to
the standard one on $\mathbb{R}^{3}$: 
\begin{equation}
\left( \mathop{\rm curl}\nolimits X\right) ^{a}=\left( \nabla
_{b}X_{c}\right) \varphi ^{abc}.  \label{curlX}
\end{equation}%
This curl operator can then also be extended to $2$-tensor $\beta $: 
\begin{equation}
\left( \mathop{\rm curl}\nolimits\beta \right) _{ab}=\left( \nabla
_{m}^{{}}\beta _{na}^{{}}\right) \varphi _{b}^{\ mn}.  \label{curlw}
\end{equation}%
Note that when $\beta _{ab}$ is symmetric, $\mathop{\rm curl}\nolimits\beta $
is traceless. It is also not difficult to see that schematically, 
\begin{equation}
\mathop{\rm curl}\nolimits\left( \left( \mathop{\rm curl}\nolimits\beta
\right) ^{t}\right) =-\Delta \beta ^{t}+\nabla \left( \mathop{\rm div}%
\nolimits\beta \right) +\mathop{\rm Riem}\nolimits\circledast \beta
+T\circledast \nabla \beta +\left( \nabla T\right) \circledast \beta
+T\circledast T\circledast \beta  \label{curlcurl}
\end{equation}%
where $\ ^{t}$ denotes transpose and $\circledast $ is some multilinear
operator involving $g,\varphi ,\psi $. From the context it will be clear
whether the curl operator is applied to a vector or a $2$-tensor.

As in \cite{GrigorianCoflow}, we can also use the $G_{2}$-structure $3$-form
to define a commutative product $\alpha \circ \beta $ of two $2$-tensors $%
\alpha $ and $\beta $%
\begin{equation}
\left( \alpha \circ \beta \right) _{ab}=\varphi _{amn}\varphi _{bpq}\alpha
^{mp}\beta ^{nq}  \label{phphprod}
\end{equation}%
Note that $\left( \alpha \circ \beta \right) ^{t}=\left( \alpha ^{t}\circ
\beta ^{t}\right) $. If $\alpha $ and $\beta $ are both symmetric or both
skew-symmetric, then $\alpha \circ \beta $ is a symmetric $2$-tensor. Also,
for a $2$-tensor we have the standard matrix product $\left( \alpha \beta
\right) _{ab}=\alpha _{a}^{\ k}\beta _{kb}$.

From \cite{CleytonIvanovClosed,GrigorianG2Torsion1,karigiannis-2007} we know
that the torsion of a $G_{2}$-structure satisfies the following
integrability condition:%
\begin{equation}
\frac{1}{2}\mathop{\rm Riem}\nolimits _{ij}^{\ \beta \gamma }\varphi _{\ \
\beta \gamma }^{\alpha }=\nabla _{i}T_{j}^{\ \alpha }-\nabla _{j}T_{i}^{\
\alpha }+T_{i}^{\ \beta }T_{j}^{\ \gamma }\varphi _{\ \beta \gamma }^{\alpha
}.  \label{dtinteg}
\end{equation}%
Taking projections of (\ref{dtinteg}) to different representations of $G_{2}$%
, we obtain the following expressions:

\begin{lemma}
\label{propTors}The torsion tensor $T$ satisfies the following identities 
\begin{subequations}%
\label{Tids}%
\begin{eqnarray}
\left( \nabla T\right) \lrcorner \psi &=&-\left( T\lrcorner \varphi \right)
\lrcorner T+T^{2}\lrcorner \varphi +\left( \mathop{\rm Tr}\nolimits T\right)
\left( T\lrcorner \varphi \right)  \label{dtabcond1} \\
0 &=&d\left( \mathop{\rm Tr}\nolimits T\right) -\mathop{\rm div}\nolimits
\left( T^{t}\right) -\left( T\lrcorner \varphi \right) \lrcorner T^{t}
\label{dtabcond2} \\
\mathop{\rm Ric}\nolimits &=&-\mathop{\rm Sym}\nolimits \left( \mathop{\rm
curl}\nolimits T^{t}-\nabla \left( T\lrcorner \varphi \right) +T^{2}-%
\mathop{\rm Tr}\nolimits \left( T\right) T\right)  \label{RicCond} \\
\frac{1}{4}\mathop{\rm Ric}\nolimits ^{\ast } &=&\mathop{\rm curl}\nolimits
T+\frac{1}{2}T\circ T  \label{RicSCond} \\
\mathrm{R} &=&2\mathop{\rm Tr}\nolimits \left( \mathop{\rm curl}\nolimits
T\right) -\psi \left( T,T\right) -\mathop{\rm Tr}\nolimits \left(
T^{2}\right) +\left( \mathop{\rm Tr}\nolimits T\right) ^{2}  \label{ScalCurv}
\end{eqnarray}%
\end{subequations}%
where $\left( \mathop{\rm Ric}\nolimits ^{\ast }\right) _{ab}=%
\mathop{\rm
Riem}\nolimits _{mnpq}^{{}}\varphi _{\ \ \ a}^{mn}\varphi _{\ \ b}^{pq}$ and 
$\psi \left( T,T\right) =\psi _{abcd}T^{ab}T^{cd}$. Note that from (\ref%
{Tdecomp}), $\mathop{\rm Tr}\nolimits T=\frac{7}{4}\tau _{0}$ and $%
T\lrcorner \varphi =-6\tau _{1}$.
\end{lemma}

The symmetric $2$-tensor $\mathop{\rm Ric}\nolimits ^{\ast }$ has been
defined and studied by Cleyton and Ivanov in \cite%
{CleytonIvanovClosed,CleytonIvanovCurv}. Note that $\mathop{\rm Tr}\nolimits
\left( \mathop{\rm Ric}\nolimits ^{\ast }\right) =2\mathrm{R} $, where $%
\mathrm{R} $ is the scalar curvature. Thus the tensors $\mathop{\rm Ric}%
\nolimits $ and $\mathop{\rm Ric}\nolimits ^{\ast }$ span the components of $%
\mathop{\rm Riem}\nolimits $ that lie in $1\oplus 27\oplus 27$
representations of $G_{2}$. It is know that $\mathop{\rm Riem}\nolimits $
has no components in the $7$ or $14$ dimensional representations of $G_{2}$.
The identities (\ref{dtabcond1}), (\ref{dtabcond2}), as well as the
projection of (\ref{RicSCond}) to $\Lambda _{14}^{2}$ are a consequence of
this. In fact, taking the skew-symmetric part of (\ref{RicSCond}) and using
the fact that $\mathop{\rm Ric}\nolimits ^{\ast }$ is by definition
symmetric, gives us 
\begin{equation}
\mathop{\rm Skew}\nolimits \left( \mathop{\rm curl}\nolimits T\right) =-%
\frac{1}{2}\mathop{\rm Skew}\nolimits \left( T\circ T\right) .
\end{equation}%
In particular, this shows that $\mathop{\rm curl}\nolimits T$ is symmetric
whenever $T$ is skew-symmetric or symmetric, and in particular, if $\varphi $
is closed or co-closed.

Let us now look at the properties of $\Delta _{\varphi }\varphi =dd^{\ast
}\varphi +d^{\ast }d\varphi $.

\begin{proposition}[\protect\cite{GrigorianCoflow}]
\label{propLap}Suppose $\varphi $ defines a $G_{2}$-structure. Then $\Delta
_{\varphi }\varphi =X\lrcorner \psi +3\mathrm{i}_{\varphi }\left( h\right) $
with%
\begin{subequations}%
\label{lapdecom}%
\begin{eqnarray}
X &=&-\mathop{\rm div}\nolimits T  \label{pi7lapphi} \\
h &=&-\frac{1}{4}\mathop{\rm Ric}\nolimits ^{\ast }+\frac{1}{6}\left( 
\mathrm{R} +2\left\vert T\right\vert ^{2}\right) g-T^{t}T-\frac{1}{2}\left(
T\lrcorner \varphi \right) \left( T\lrcorner \varphi \right)
\label{pi27lapphi} \\
&&+\frac{1}{4}T\circ T+\frac{1}{4}T^{t}\circ T^{t}-\frac{1}{2}T\circ T^{t}+%
\mathop{\rm Sym}\nolimits \left( \left( T\right) \left( T\lrcorner \psi
\right) -\left( T^{t}\right) \left( T\lrcorner \psi \right) \right) .  \notag
\end{eqnarray}%
\end{subequations}%
In particular, 
\begin{equation}
\mathop{\rm Tr}\nolimits h=\frac{2}{3}\mathrm{R} +\frac{4}{3}\left\vert
T\right\vert ^{2}.  \label{Trh}
\end{equation}
\end{proposition}

The leading order terms in $\Delta _{\varphi }\varphi $ are those that
contain second derivatives of $\varphi $, and hence first derivatives of $T$%
. Thus, $\mathop{\rm div}\nolimits T$ fully defines the $\Lambda _{7}^{3}$
component of $\Delta _{\varphi }\varphi $ and the leading order terms in $%
\Lambda _{1\oplus 27}^{3}$ are given by 
\begin{equation}
-\frac{1}{4}\mathop{\rm Ric}\nolimits ^{\ast }+\frac{1}{6}\mathop{\rm R}%
\nolimits g\sim -\mathop{\rm curl}\nolimits T+\frac{1}{3}\mathop{\rm Tr}%
\nolimits \left( \mathop{\rm curl}\nolimits T\right) g.  \label{lapphi27lo}
\end{equation}

\section{Flows of $G_{2}$-structures}

\label{SecDeformPhi}\setcounter{equation}{0}Suppose $\varphi \left( t\right) 
$ is a one-parameter family of $G_{2}$-structures on a manifold $M$ that
satisfies 
\begin{equation}
\frac{\partial \varphi \left( t\right) }{\partial t}=X\left( t\right)
\lrcorner \psi \left( t\right) +3\mathrm{i}_{\varphi \left( t\right) }\left(
h\left( t\right) \right) .  \label{dphit}
\end{equation}%
As shown by Karigiannis in \cite{karigiannis-2007}, the associated
quantities $g\left( t\right) ,\mathop{\rm vol}\nolimits _{t},\psi \left(
t\right) $, $T\left( t\right) $ satisfy the following evolution equations:

\begin{lemma}[\protect\cite{karigiannis-2007}]
If $\varphi \left( t\right) $ satisfies the equation (\ref{dphit}), then we
also have the following equations: 
\begin{subequations}%
\label{evoleqs}%
\begin{eqnarray}
\frac{\partial g}{\partial t} &=&2h  \label{gevol} \\
\frac{\partial \mathop{\rm vol}\nolimits }{\partial t} &=&\mathop{\rm Tr}%
\nolimits \left( h\right) \mathop{\rm vol}\nolimits  \label{volevol} \\
\frac{\partial \psi }{\partial t} &=&4\mathrm{i}_{\psi }\left( h\right)
-X\wedge \varphi  \label{psievol} \\
\frac{\partial T}{\partial t} &=&\nabla X-\mathop{\rm curl}\nolimits
h+Th-\left( T\right) \left( X\lrcorner \varphi \right)  \label{Tevol}
\end{eqnarray}%
\end{subequations}%
where $\mathrm{i}_{\psi }\left( h\right) _{abcd}=-h_{[a}^{e}\psi
_{bcd]e}^{{}}\ $and equivalently, $4\mathrm{i}_{\psi }\left( h\right)
=-3\ast \mathrm{i}_{\varphi }\left( h\right) +\left( \mathop{\rm Tr}%
\nolimits h\right) \psi $.
\end{lemma}

Similarly, as in \cite{GrigorianCoflow}, we can consider flows of $\psi ,$
given by 
\begin{equation}
\frac{\partial \psi \left( t\right) }{\partial t}=\ast \left( X\left(
t\right) \lrcorner \psi \left( t\right) \right) +3\ast \mathrm{i}_{\varphi
\left( t\right) }\left( s\left( t\right) \right)  \label{psiflow}
\end{equation}%
for some symmetric $2$-tensor $s$. Since $3\ast \mathrm{i}_{\varphi }\left(
s\right) =4\mathrm{i}_{\psi }\left( \frac{1}{4}\left( \mathop{\rm Tr}%
\nolimits s\right) g-s\right) $, comparing (\ref{psiflow}) with (\ref%
{psievol}) give us corresponding evolution equations for $\varphi \left(
t\right) $, $g\left( t\right) ,\mathop{\rm vol}\nolimits _{t},T\left(
t\right) $ from (\ref{dphit}) and (\ref{evoleqs}) by taking $h=\frac{1}{4}%
\left( \mathop{\rm Tr}\nolimits s\right) g-s$.

When constructing geometric flows, there are two main considerations: 1) the
flow's stationary points should correspond to geometrically interesting
objects; and 2) the flow should be parabolic in some sense. The first
property is the main motivation for studying a flow, since we ideally want
the flow to deform a geometric structure to one that has nicer or more
constrained properties and the second property is a minimal requirement to
at least guarantee short-time existence and uniqueness of solutions. In \cite%
{GaoChen1}, Chen defined a class of \emph{reasonable }flows (\ref{dphit}) of 
$G_{2}$-structures that satisfy the following $4$ general conditions:

\begin{enumerate}
\item The metric should evolve by the Ricci flow to leading order, and be no
more than quadratic in the torsion, that is 
\begin{equation}
\frac{\partial g}{\partial t}=2h=-2\mathop{\rm Ric}\nolimits +Cg+L\left(
T\right) +T\circledast T  \label{dgtreason}
\end{equation}%
where $C$ is a constant and $L$ is some linear operator involving $g,\varphi
,\psi $.

\item The vector field $X$ is at most linear in $\nabla T$ and at most
quadratic in $T$:%
\begin{equation}
X=L\left( \nabla T\right) +L\left( T\right) +L\left( \mathop{\rm Riem}%
\nolimits \right) +T\circledast T+C.  \label{Xgen}
\end{equation}

\item The torsion tensor should evolve by $\Delta T$ to leading order, and
be at most linear in $\mathop{\rm Riem}\nolimits $ and $\nabla T,$ and at
most cubic in $T$:%
\begin{eqnarray}
\frac{\partial T}{\partial t} &=&\Delta T+L\left( \nabla T\right) +L\left( %
\mathop{\rm Riem}\nolimits \right) +\mathop{\rm Riem}\nolimits \circledast
T+\nabla T\circledast T  \label{dTdt} \\
&&+L\left( T\right) +T\circledast T+T\circledast T\circledast T.  \notag
\end{eqnarray}

\item The flow (\ref{dphit}) has short-time existence and uniqueness.
\end{enumerate}

As one of the key properties of \emph{reasonable }flows defined above is
that the flow of the metric is the Ricci flow to leading order, we will
instead refer to flows that satisfy properties 1.-4. as \emph{Ricci-like
flows}. This is appropriate because a variety of techniques that originated
from the study of the Ricci flow have been applied to these flows. In
particular, under the Ricci flow, invariants of the metric $\mathop{\rm Riem}%
\nolimits,$ $\mathop{\rm Ric}\nolimits,$ $\mathrm{R},$ all satisfy heat-like
equations. Therefore it is appropriate that for a Ricci-like flow of a $G_{2}
$-structure, the torsion, which an invariant of the $G_{2}$-structure also
satisfies a heat-like equation (\ref{dTdt}). This is important because then $%
\nabla ^{k}T$ and $\left\vert T\right\vert ^{2}$ also satisfy heat-like
equations and this is necessary to be able to obtain estimates using the
maximum principle.

Using techniques developed by Shi in \cite{ShiEstimate} for the Ricci flow
and their adaptation to $G_{2}$-structures by Lotay and Wei \cite{LotayWei1}%
, Chen then showed that a reasonable flow satisfies the following Shi-type
estimate.

\begin{theorem}[{\protect\cite[Theorem 2.1]{GaoChen1}}]
\label{ThmShi}Suppose (\ref{dphit}) is a Ricci-like flow of $G_{2}$%
-structures, such that the coefficients in equations (\ref{dphit}), (\ref%
{dgtreason}), (\ref{Xgen}), and (\ref{dTdt}) are bounded by a constant $%
\Lambda $. Let $B_{r}\left( p\right) $ be a ball of radius $r$ with respect
to the initial metric $g\left( 0\right) $. If 
\begin{equation}
\left\vert \mathop{\rm Riem}\nolimits \left( x,t\right) \right\vert
_{g\left( t\right) }+\left\vert T\left( x,t\right) \right\vert _{g\left(
t\right) }^{2}+\left\vert \nabla T\left( x,t\right) \right\vert _{g\left(
t\right) }<\Lambda  \label{ShiEst1}
\end{equation}%
for any $\left( x,t\right) \in $ $B_{r}\left( p\right) \times \left[ 0,t_{0}%
\right] $, then 
\begin{equation}
\left\vert \nabla ^{k}\mathop{\rm Riem}\nolimits \left( x,t\right)
\right\vert _{g\left( t\right) }+\left\vert \nabla ^{k+1}T\left( x,t\right)
\right\vert _{g\left( t\right) }<C\left( k,r,\Lambda ,t\right)
\label{ShiEst2}
\end{equation}%
for any $\left( x,t\right) \in $ $B_{r/2}\left( p\right) \times \left[ \frac{%
t_{0}}{2},t_{0}\right] $ for all $k=1,2,3,...$
\end{theorem}

It should be noted that in \cite{LotayWei1}, the condition analogous to (\ref%
{ShiEst1}) does not include a $\left\vert T\right\vert ^{2}$ term. This is
because in the case of a closed $G_{2}$-structure, $\left\vert T\right\vert
^{2}=-\mathrm{R} \leq C\left\vert \mathop{\rm Riem}\nolimits \right\vert .$
Therefore, the norm of the torsion can always be bounded in terms of the
norm of $\mathop{\rm Riem}\nolimits .$ For other torsion classes, and in
particular, co-closed $G_{2}$-structures, this is no longer true, therefore $%
\left\vert T\right\vert ^{2}$ needs to be included in (\ref{ShiEst1}).

Using the estimates from Theorem \ref{ThmShi}, Chen then derived an estimate
for the blow-up rate on a compact manifold.

\begin{theorem}[{\protect\cite[Theorem 5.1]{GaoChen1}}]
\label{ThmBlowUp}If $\varphi \left( t\right) $ is a solution to a Ricci-like
flow of $G_{2}$-structures on a compact manifold in a finite maximal time
interval $[0,t_{0})$, then 
\begin{equation}
\sup_{M}\left( \left\vert \mathop{\rm Riem}\nolimits \left( x,t\right)
\right\vert _{g\left( t\right) }^{2}+\left\vert T\left( x,t\right)
\right\vert _{g\left( t\right) }^{4}+\left\vert \nabla T\left( x,t\right)
\right\vert _{g\left( t\right) }^{2}\right) ^{\frac{1}{2}}\geq \frac{C}{%
t_{0}-t}  \label{blowupestimate}
\end{equation}%
for some positive constant $C$.
\end{theorem}

The estimate (\ref{blowupestimate}) shows that a solution will exist as long
the quantity of the left-hand side of (\ref{blowupestimate}) remains bounded.

A classic example of a Ricci-like flow of $G_{2}$-structures is the
Laplacian flow of $G_{2}$-structures that was introduced by Bryant in \cite%
{bryant-2003}: 
\begin{equation}
\frac{\partial \varphi }{\partial t}=\Delta _{\varphi }\varphi .
\label{lapflow}
\end{equation}%
If the initial $G_{2}$-structure is closed, then this property is preserved
along the flow. It is then natural to think of (\ref{lapflow}) as a flow of
closed $G_{2}$-structures. In this case, since $T^{t}=-T$, from (\ref{Tids}%
), $\mathop{\rm Ric}\nolimits^{\ast }=4\mathop{\rm Ric}\nolimits+T%
\circledast T$ and $\mathrm{R}=2\mathop{\rm Tr}\nolimits\left( \mathop{\rm
curl}\nolimits T\right) -\psi \left( T,T\right) -\mathop{\rm Tr}%
\nolimits\left( T^{2}\right) =-\left\vert T\right\vert ^{2}$; and thus, from
(\ref{pi27lapphi})$,$ $h=-\mathop{\rm Ric}\nolimits+T\circledast T$, and so
from (\ref{gevol}), we do find that (\ref{dgtreason}) holds. Moreover, from (%
\ref{dtabcond2}), we see that $\mathop{\rm div}\nolimits T=0$ in this case,
and hence $X=0$. The expression (\ref{dTdt}) comes from (\ref{Tevol}) and
using $h=-\mathop{\rm curl}\nolimits T+T\circledast T$ 
\begin{equation}
\frac{\partial T}{\partial t}=\mathop{\rm curl}\nolimits\left( \mathop{\rm
curl}\nolimits T\right) +\nabla T\circledast T+T\circledast T\circledast T.
\label{dTclosed}
\end{equation}%
Using (\ref{curlcurl}) to expand $\mathop{\rm curl}\nolimits\left( %
\mathop{\rm curl}\nolimits T\right) $ together the facts that $%
\mathop{\rm
curl}\nolimits T$ is symmetric, $T$ is skew-symmetric, and $\mathop{\rm div}%
\nolimits T=0$, allows to express the right-hand side of (\ref{dTclosed}) as 
$\Delta T+\mathop{\rm Riem}\nolimits\circledast T+\nabla T\circledast
T+T\circledast T\circledast T$. Finally, short-term existence and uniqueness
of the flow (\ref{lapflow}) has been first proved by Bryant and Xu in \cite%
{BryantXu}. For more on the properties of this flow, as well as the details
of the above calculations, the reader is referred to the series of papers by
Lotay and Wei \cite{LotayWei1,LotayWei2,LotayWei1a}. The results in Theorems %
\ref{ThmShi} and \ref{ThmBlowUp} are extensions of similar results for the
Laplacian flow of closed $G_{2}$-structures in \cite{LotayWei1}.

\section{Laplacian coflow}

\label{secLapcoflow}\setcounter{equation}{0}In \cite{KarigiannisMcKayTsui},
Karigiannis, McKay, and Tsui introduced an alternative flow of $G_{2}$%
-structures, called the Laplacian \emph{coflow}:%
\begin{equation}
\frac{\partial \psi }{\partial t}=-\Delta _{\psi }\psi .  \label{cofloworig}
\end{equation}%
If the initial $G_{2}$-structure is co-closed, i.e. $d\psi =0$, then this
property is preserved along the flow. Therefore, the coflow may be regarded
as a natural flow of co-closed $G_{2}$-structures. In order to understand
flows of co-closed $G_{2}$-structures, we need to understand better the
properties of $T$ and the Hodge Laplacian in this case. Rewriting Lemma \ref%
{propTors} and Proposition \ref{propLap} in the case of symmetric $T$, we
find the following.

\begin{proposition}
Suppose $\varphi $ is a co-closed $G_{2}$-structure, then the torsion tensor 
$T$ satisfies the following identities 
\begin{subequations}%
\label{TidsCocl}%
\begin{eqnarray}
\mathop{\rm div}\nolimits T &=&d\left( \mathop{\rm Tr}\nolimits T\right) \\
\mathop{\rm curl}\nolimits T &=&\left( \mathop{\rm curl}\nolimits T\right)
^{t} \\
\mathop{\rm Ric}\nolimits &=&\mathop{\rm curl}\nolimits T-T^{2}+\mathop{\rm
Tr}\nolimits \left( T\right) T \\
\frac{1}{4}\mathop{\rm Ric}\nolimits ^{\ast } &=&\mathop{\rm curl}\nolimits
T+\frac{1}{2}T\circ T=\mathop{\rm Ric}\nolimits +\frac{1}{2}T\circ T+T^{2}-%
\mathop{\rm Tr}\nolimits \left( T\right) T \\
\mathrm{R} &=&\left( \mathop{\rm Tr}\nolimits T\right) ^{2}-\left\vert
T\right\vert ^{2}.
\end{eqnarray}%
\end{subequations}%
The Hodge Laplacian is given by $\Delta _{\varphi }\varphi =X\lrcorner \psi
+3\mathrm{i}_{\varphi }\left( s\right) $ with%
\begin{subequations}%
\label{lapdecomcocl}%
\begin{eqnarray}
X &=&-\mathop{\rm div}\nolimits T  \label{Xcoclosed} \\
s &=&-\mathop{\rm Ric}\nolimits +\frac{1}{6}\left( \mathrm{R} +2\left\vert
T\right\vert ^{2}\right) g+\mathop{\rm Tr}\nolimits \left( T\right) T-2T^{2}-%
\frac{1}{2}T\circ T \\
&=&-\mathop{\rm curl}\nolimits T+\frac{1}{6}\left( \left( \mathop{\rm Tr}%
\nolimits T\right) ^{2}+\left\vert T\right\vert ^{2}\right) g-T^{2}-\frac{1}{%
2}T\circ T \\
\mathop{\rm Tr}\nolimits s &=&\frac{2}{3}\mathrm{R} +\frac{4}{3}\left\vert
T\right\vert ^{2}=\frac{2}{3}\left( \left( \mathop{\rm Tr}\nolimits T\right)
^{2}+\left\vert T\right\vert ^{2}\right) .
\end{eqnarray}%
\end{subequations}%
\end{proposition}

Comparing (\ref{cofloworig}) with (\ref{psiflow}) and using (\ref%
{lapdecomcocl}), we see that to leading order the evolution of the metric is
given by $2\mathop{\rm Ric}\nolimits,$ that is the opposite of the Ricci
flow. Thus, in order for the flow to be Ricci-like and to have any hope of
existence and uniqueness, the sign in (\ref{cofloworig}) needs to be
reversed. Therefore, let us redefine the Laplacian coflow as 
\begin{equation}
\frac{d\psi }{dt}=\Delta _{\psi }\psi .  \label{lapcoflow1}
\end{equation}%
We then find that 
\begin{equation}
\frac{\partial g}{\partial t}=-2\mathop{\rm Ric}\nolimits+T\circ T+2\left( %
\mathop{\rm Tr}\nolimits T\right) T  \label{coclosedgdt2}
\end{equation}%
which now satisfies (\ref{dgtreason}). Also, $X=-\mathop{\rm div}\nolimits T$%
, which satisfies (\ref{Xgen}). To obtain the general form of the evolution
of the torsion, note that to leading order, $h=-s=\mathop{\rm curl}\nolimits
T$, so from (\ref{Tevol}), 
\begin{equation*}
\frac{\partial T}{\partial t}=-\nabla \left( \mathop{\rm div}\nolimits
T\right) -\mathop{\rm curl}\nolimits\left( \mathop{\rm curl}\nolimits
T\right) +\nabla T\circledast T
\end{equation*}%
however, since both $T$ and $\mathop{\rm curl}\nolimits T$ are symmetric, 
\begin{equation*}
\mathop{\rm curl}\nolimits\left( \mathop{\rm curl}\nolimits T\right)
=-\Delta T+\nabla \left( \mathop{\rm div}\nolimits T\right) +\mathop{\rm
Riem}\nolimits\circledast T+\left( \nabla T\right) \circledast
T+T\circledast T\circledast T
\end{equation*}%
Hence, overall, 
\begin{equation}
\frac{\partial T}{\partial t}=\Delta T-2\nabla \left( \mathop{\rm div}%
\nolimits T\right) +\mathop{\rm Riem}\nolimits\circledast T+\left( \nabla
T\right) \circledast T+T\circledast T\circledast T.  \label{DTdtcoflow}
\end{equation}%
Notice that this does not satisfy (\ref{dTdt}). In fact, we can see that the
presence of the $\nabla \left( \mathop{\rm div}\nolimits T\right) $ term in (%
\ref{DTdtcoflow}) is due to the negative sign of $\mathop{\rm div}\nolimits
T $ in (\ref{Xcoclosed}). As it was shown in \cite{GrigorianCoflow}, the
sign of $\mathop{\rm div}\nolimits T$ also causes problems at a much more
fundamental level: it prevents the flow (\ref{lapcoflow1}) from being
parabolic even along closed $4$-forms. Proposition \ref{propLappsilin} below
gives the linearization of $\Delta _{\psi }$. It is then easy to see that
for closed $4$-forms, the symbol will be negative in the $\Lambda _{7}^{4}$
direction, but non-negative in $\Lambda _{27}^{4}.$

\begin{proposition}[{\protect\cite[Prop. 4.7]{GrigorianCoflow}}]
\label{propLappsilin}The linearization of $\Delta _{\psi }$ at $\psi $ is
given by 
\begin{subequations}%
\label{laplincocl} 
\begin{eqnarray}
\pi _{7}\left( D_{\psi }\Delta _{\psi }\right) \left( \chi \right)
&=&d\left( \mathop{\rm div}\nolimits X\right) \wedge \varphi +l.o.t. \\
\pi _{1\oplus 27}\left( D_{\psi }\Delta _{\psi }\right) \left( \chi \right)
&=&\frac{3}{2}\ast \mathrm{i}_{\varphi }\left( \Delta h+\frac{1}{4}%
\mathop{\rm Hess}\nolimits \left( \mathop{\rm Tr}\nolimits h\right) -\frac{1%
}{2}\left( \Delta \mathop{\rm Tr}\nolimits h\right) g\right. \\
&&\left. -\mathop{\rm Sym}\nolimits \left( \nabla \mathop{\rm div}\nolimits
h+\mathop{\rm curl}\nolimits \left( \nabla X\right) ^{t}\right)
+l.o.t.\right)  \notag
\end{eqnarray}%
\end{subequations}%
where $\chi =\ast \left( X\lrcorner \psi +3\mathrm{i}_{\varphi }\left(
h\right) \right) $. Moreover, if $\chi $ is closed, we can write $D_{\psi
}\Delta _{\psi }$ as 
\begin{equation}
D_{\psi }\Delta _{\psi }\left( \chi \right) =-\Delta _{\psi }\chi -\mathcal{L%
}_{V\left( \chi \right) }\psi +2d\left( \left( \mathop{\rm div}\nolimits
X\right) \varphi \right) +dF\left( \chi \right)  \label{coflowlin}
\end{equation}%
where 
\begin{equation}
V\left( \chi \right) =\frac{3}{4}\nabla \mathop{\rm Tr}\nolimits h-2%
\mathop{\rm curl}\nolimits X  \label{coflowvect}
\end{equation}%
and $F\left( \chi \right) $ is a $3$-form-valued algebraic function of $\chi 
$.
\end{proposition}

Looking closer at the leading terms in the linearization (\ref{coflowlin})
evaluated at closed forms, we see that the term $2d\left( \left( 
\mathop{\rm
div}\nolimits X\right) \varphi \right) $ appears for exactly the same reason
as the term $-2\nabla \left( \mathop{\rm div}\nolimits T\right) $ in (\ref%
{DTdtcoflow}) - namely the \textquotedblleft wrong\textquotedblright\ sign
of the $\pi _{7}$ component of $\Delta _{\psi }\psi $. To fix this problem,
in \cite{GrigorianCoflow}, a \emph{modified Laplacian coflow} has been
proposed: 
\begin{equation}
\frac{\partial \psi }{\partial t}=\Delta _{\psi }\psi +2d\left( \left( A-%
\mathop{\rm Tr}\nolimits T\right) \varphi \right)  \label{coflowmod2}
\end{equation}%
where $A$ is some constant. Since for co-closed $G_{2}$-structures, $%
\mathop{\rm Tr}\nolimits T=\mathop{\rm div}\nolimits T,$ the leading term in
the modification precisely reverses the sign of the $\Lambda _{7}^{4}$
component of the original flow (\ref{lapcoflow}). However, because we want
the right hand side of the flow to be an exact $4$-form for co-closed $G_{2}$%
-structures, there are some additional lower order terms. The constant $A$
could be set to zero, however adding it may allow for more flexibility. The
linearization of the modified coflow at a closed $4$-form is now given by%
\begin{equation}
\frac{\partial \chi }{\partial t}=-\Delta _{\psi }\chi -\mathcal{L}_{V\left(
\chi \right) }\psi +d\hat{F}\left( \chi \right)  \label{modcoflowlin}
\end{equation}%
where $V\left( \chi \right) $ is as in (\ref{coflowvect}) and $\hat{F}\left(
\chi \right) $ involves no derivatives of $\chi $. Hence, in the direction
of closed forms, this flow is now weakly parabolic. Moreover, the undesired
term is removed from the evolution equation for $T$ and its evolution is now
given by (\ref{dTdt}).

The additional term in (\ref{coflowmod2}) now also allows to prove
short-time existence and uniqueness, hence completing the requirements for (%
\ref{coflowmod2}) to be a Ricci-like flow. The proof, as given in \cite%
{GrigorianCoflow}, follows a procedure similar to the approach taken by
Bryant and Xu \cite{BryantXu} for the proof of short-time existence and
uniqueness for the Laplacian flow (\ref{lapflow}), which is in turn based on
DeTurck's \cite{DeTurckTrick} and Hamilton's \cite{HamiltonNashMoser}
approaches to the proof of short-time existence and uniqueness of the Ricci
flow. Let $\psi \left( t\right) =\psi _{0}+\chi \left( t\right) $ where $%
\chi \left( t\right) $ is an exact $4$-form with $\chi \left( 0\right) =0$.
Then, given this initial condition, the flow (\ref{coflowmod2}) can be
rewritten as an initial value problem for $\chi \left( t\right) $. From the
linearization (\ref{modcoflowlin}) we see that by adding the term $\mathcal{L%
}_{V\left( \chi \left( t\right) \right) }\psi \left( t\right) $ we obtain a
strictly parabolic flow in the direction of closed forms, which is related
to the original flow by diffeomorphism:%
\begin{equation}
\frac{\partial \chi }{\partial t}=\Delta _{\psi }\psi +2d\left( \left( A-%
\mathop{\rm Tr}\nolimits T_{\psi }\right) \ast _{\psi }\psi \right) +%
\mathcal{L}_{V\left( \chi \right) }\psi .  \label{modcoflow3}
\end{equation}

This is the essence of what is known as \textquotedblleft DeTurck's
trick\textquotedblright\ - turning a weakly parabolic flow into a strictly
parabolic one. In the case of Ricci flow this is enough to obtain short-time
existence and uniqueness, however in this case, the parabolicity is only
along closed forms, hence we cannot apply the standard parabolic theory
right away, and more steps are needed. Let us also define the spaces of
time-dependent and time-independent exact $4$-forms $\mathcal{F}$ and $%
\mathcal{G}$, respectively. Moreover, since we know that $\psi \left(
t\right) $ always defines a $G_{2}$-structure and is thus a positive $4$%
-form, $\chi $ will always lie in an open subset $\mathcal{U}\subset 
\mathcal{F}$ defined by 
\begin{equation}
\mathcal{U}=\left\{ \chi \in \mathcal{F}:\psi _{0}+\chi \ \ \text{is a
positive }4\text{-form}\right\} .  \label{Uspace}
\end{equation}%
Moreover, let us now define a map $F:\mathcal{U}\longrightarrow \mathcal{F}%
\times \mathcal{G}$ given by 
\begin{equation}
\chi \longrightarrow \left( \frac{\partial \chi }{\partial t}-\Delta _{\psi
}\psi -2d\left( \left( A-\mathop{\rm Tr}\nolimits T_{\psi }\right) \ast
_{\psi }\psi \right) -\mathcal{L}_{V\left( \chi \right) }\psi ,\left. \chi
\right\vert _{t=0}\right) .  \label{mapF}
\end{equation}%
Adapting the results in \cite{BryantXu}, it is easy to see $\mathcal{F}$, $%
\mathcal{G}$, and $\mathcal{H}:=\mathcal{F}\times \mathcal{G}$ are \emph{%
graded tame Fr\'{e}chet spaces}. Moreover, it was then shown in \cite%
{GrigorianCoflow} that $F$ is smooth \emph{tame }map of Fr\'{e}chet spaces,
such that its derivative $DF\left( \chi \right) :\mathcal{F}\longrightarrow 
\mathcal{H}$ is an isomorphism for all $\chi \in \mathcal{U}$ and the
inverse $\left( DF\right) ^{-1}:\mathcal{U}\times \mathcal{H}\longrightarrow 
\mathcal{F}$ is smooth tame. The significance of these facts are that in the
category of Fr\'{e}chet spaces there exists an inverse function theorem -
the Nash-Moser Inverse Function \cite{HamiltonNashMoser}, which tells us
that the map $F$ is locally invertible. From this it follows that the flow (%
\ref{modcoflow3}) has short-time existence and uniqueness.

To prove short-time existence and uniqueness for the flow (\ref{coflowmod2})
we need to relate (\ref{coflowmod2}) and (\ref{modcoflow3}) via
diffeomorphisms. Suppose $\bar{\chi}\left( t\right) $ is the unique
short-time solution to (\ref{modcoflow3}), and $\bar{\psi}=\psi _{0}+\bar{%
\chi}$. Consider the following ODE for a family of diffeomorphisms $\phi
_{t} $:%
\begin{equation}
\left\{ 
\begin{array}{c}
\frac{\partial \phi _{t}}{\partial t}=-V\left( \bar{\chi}\left( t\right)
\right) \\ 
\phi _{0}=\mathop{\rm id}\nolimits%
\end{array}%
\right.  \label{ddtdiffeo2}
\end{equation}%
This has a unique solution $\phi _{t}$. Now let $\psi \left( t\right)
=\left( \phi _{t}\right) ^{\ast }\bar{\psi}\left( t\right) $, then $\psi
\left( 0\right) =\psi _{0}$, and since diffeomorphisms commute with $d$, $%
\psi \left( t\right) $ is closed for all $t$. Moreover, as shown in \cite[%
Theorem 6.9]{GrigorianCoflow}, $\psi \left( t\right) $ now satisfies (\ref%
{coflowmod2}). Uniqueness is obtained similarly using the uniqueness of
solutions of (\ref{ddtdiffeo2}). Hence, overall, we obtain a unique
short-time solution for the modified Laplacian coflow (\ref{coflowmod2}) and
can now conclude that it is a Ricci-like flow.

\begin{theorem}
The Laplacian coflow (\ref{coflowmod2}) of co-closed $G_{2}$-structures is a
Ricci-like flow.
\end{theorem}

\section{Further directions}

There are several important unanswered questions regarding flows of
co-closed $G_{2}$-structures. An intriguing question is whether it is
possible to obtain at least short-time existence and uniqueness of the
unmodified Laplacian coflow (\ref{lapcoflow}). To leading order the only
difference with the modified coflow is the sign of the $\Lambda _{7}^{4}$
component which is given by $\mathop{\rm div}\nolimits T$. So in particular,
if $\mathop{\rm div}\nolimits T$ vanishes, then the two flows agree. It is
also known \cite{karigiannis-2005-57} that deformations in the $\Lambda
_{7}^{4}$ directions keep the metric unchanged. Moreover, in \cite%
{GrigorianOctobundle}, the torsion $T$ has been shown to play a role of an
octonionic connection on the bundle of $G_{2}$-structures that correspond to
the same metric, which can be given the structure of an octonion bundle. In
this interpretation, on a compact manifold, the condition $\mathop{\rm div}%
\nolimits T=0$ corresponds to critical points of the functional $\int
\left\vert T\right\vert ^{2}\mathop{\rm vol}\nolimits ,$ and is hence the
analog of a Coulomb gauge. It is therefore tempting to think that to relate
the flows (\ref{lapcoflow}) and (\ref{mod-coflow}), a gauge-fixing condition
such as $\mathop{\rm div}\nolimits T=0$ needs to be introduced.

There are also multiple questions relating to the modified coflow itself. As
it is a Ricci-like flow, Shi-type estimates apply to it, so it is likely
that in addition to Chen's results in \cite{GaoChen1}, more properties such
as real analyticity and stability could be proved using techniques similar
to the ones used by Lotay and Wei in \cite{LotayWei1,LotayWei2,LotayWei1a}.
Indeed, as this article was being finalized, the author was made aware that
Bedulli and Vezzoni \cite{BedulliVezzoni} have generalized the proof of
stability from \cite{LotayWei1a} to a wider class of geometric flows that
also includes the modified Laplacian coflow with $A=0$.

Apart from the Laplacian flow and the coflows, there could be more
interesting flows of $G_{2}$-structures. For co-closed $G_{2}$-structures,
it is an open question whether the flow $\frac{\partial \varphi }{\partial t}%
=d^{\ast }d\varphi $ satisfies the co-closed condition. More generally, the
conditions for a flow to be Ricci-like is a good set of conditions that
flows should satisfy. In particular, one could try to construct flows using
the first 3 conditions, but then also making sure that short-time existence
and uniqueness is satisfied.

\bibliographystyle{jhep-a}
\bibliography{refs2}

\providecommand{\href}[2]{#2}\begingroup\raggedright\begin{thebibliography}{10}

\bibitem{BagagliniFernandezFino1}
L.~Bagaglini, M.~Fern{\'a}ndez and A.~Fino, {\it Laplacian coflow on the
  7-dimensional {H}eisenberg group},
  \href{http://arXiv.org/abs/1704.00295}{{\tt 1704.00295}}.

\bibitem{BagagliniFino1}
L.~Bagaglini and A.~Fino, {\it The laplacian coflow on almost-abelian {L}ie
  groups},  \href{http://arXiv.org/abs/1711.03751}{{\tt 1711.03751}}.

\bibitem{BedulliVezzoni}
L.~Bedulli and L.~Vezzoni, {\it Stability of geometric flows of closed forms},
  \href{http://arXiv.org/abs/1811.09416}{{\tt 1811.09416}}.

\bibitem{Bryant-1987}
R.~L. Bryant, {\it Metrics with exceptional holonomy},  {\em Ann. of Math. (2)}
  {\bf 126} (1987), no.~3 525--576.

\bibitem{bryant-2003}
R.~L. Bryant, {\it Some remarks on {G}\_2-structures},  in {\em Proceedings of
  {G}\"okova {G}eometry-{T}opology {C}onference 2005}, pp.~75--109, G\"okova
  Geometry/Topology Conference (GGT), G\"okova, 2006.
\newblock \href{http://arXiv.org/abs/math/0305124}{{\tt math/0305124}}.

\bibitem{BryantXu}
R.~L. Bryant and F.~Xu, {\it {Laplacian Flow for Closed $G_2$-Structures: Short
  Time Behavior}},  \href{http://arXiv.org/abs/1101.2004}{{\tt 1101.2004}}.

\bibitem{GaoChen1}
G.~Chen, {\it Shi-type estimates and finite time singularities of flows of
  {$G_2$} structures},  {\em Q. J. Math.} (2018)
  [\href{http://arXiv.org/abs/1703.08526}{{\tt 1703.08526}}].

\bibitem{CleytonIvanovClosed}
R.~Cleyton and S.~Ivanov, {\it On the geometry of closed {$G_2$}-structures},
  {\em Comm. Math. Phys.} {\bf 270} (2007), no.~1 53--67
  [\href{http://arXiv.org/abs/math/0306362}{{\tt math/0306362}}].

\bibitem{CleytonIvanovCurv}
R.~Cleyton and S.~Ivanov, {\it Curvature decomposition of {$G_2$}-manifolds},
  {\em J. Geom. Phys.} {\bf 58} (2008), no.~10 1429--1449.

\bibitem{CrowleyNordstrom1}
D.~Crowley and J.~Nordstr\"om, {\it New invariants of {$G_2$}-structures},
  {\em Geom. Topol.} {\bf 19} (2015), no.~5 2949--2992
  [\href{http://arXiv.org/abs/1211.0269}{{\tt 1211.0269}}].

\bibitem{DeTurckTrick}
D.~M. DeTurck, {\it Deforming metrics in the direction of their {R}icci
  tensors},  {\em J. Differential Geom.} {\bf 18} (1983), no.~1 157--162.

\bibitem{FernandezGray}
M.~Fern{\'a}ndez and A.~Gray, {\it Riemannian manifolds with structure group
  {$G\sb{2}$}},  {\em Ann. Mat. Pura Appl. (4)} {\bf 132} (1982) 19--45.

\bibitem{GrigorianG2TorsionWarped}
S.~Grigorian, {\it {G2-structure deformations and warped products}},  in {\em
  String-Math 2011}, Proceedings of Symposia in Pure Mathematics, AMS, 2012.
\newblock \href{http://arXiv.org/abs/1110.4594}{{\tt 1110.4594}}.

\bibitem{GrigorianCoflow}
S.~Grigorian, {\it Short-time behaviour of a modified {L}aplacian coflow of
  {G2}-structures},  {\em Adv. Math.} {\bf 248} (2013) 378--415
  [\href{http://arXiv.org/abs/1209.4347}{{\tt 1209.4347}}].

\bibitem{GrigorianG2Torsion1}
S.~Grigorian, {\it Deformations of {$G_2$}-structures with torsion},  {\em
  Asian J. Math.} {\bf 20} (2016), no.~1 123--155
  [\href{http://arXiv.org/abs/1108.2465}{{\tt 1108.2465}}].

\bibitem{GrigorianSU3flow}
S.~Grigorian, {\it Modified {L}aplacian coflow of {$G_2$}-structures on
  manifolds with symmetry},  {\em Differential Geom. Appl.} {\bf 46} (2016)
  39--78 [\href{http://arXiv.org/abs/1504.05506}{{\tt 1504.05506}}].

\bibitem{GrigorianOctobundle}
S.~Grigorian, {\it ${G}_2$-structures and octonion bundles},  {\em Adv. Math.}
  {\bf 308} (2017) 142--207 [\href{http://arXiv.org/abs/1510.04226}{{\tt
  1510.04226}}].

\bibitem{GrigorianYau1}
S.~Grigorian and S.-T. Yau, {\it {Local geometry of the {G}2 moduli space}},
  {\em Comm. Math. Phys.} {\bf 287} (2009) 459--488
  [\href{http://arXiv.org/abs/0802.0723}{{\tt 0802.0723}}].

\bibitem{HamiltonNashMoser}
R.~S. Hamilton, {\it The inverse function theorem of {N}ash and {M}oser},  {\em
  Bull. Amer. Math. Soc. (N.S.)} {\bf 7} (1982), no.~1 65--222.

\bibitem{Hamilton3folds}
R.~S. Hamilton, {\it Three-manifolds with positive {R}icci curvature},  {\em J.
  Differential Geom.} {\bf 17} (1982), no.~2 255--306.

\bibitem{Hitchin:2000jd-arxiv}
N.~J. Hitchin, {\it The geometry of three-forms in six and seven dimensions},
  {\em J. Differential Geom.} {\bf 55} (2000), no.~3 547--576
  [\href{http://arXiv.org/abs/math/0010054}{{\tt math/0010054}}].

\bibitem{Joycebook}
D.~D. Joyce, {\em Compact manifolds with special holonomy}.
\newblock Oxford Mathematical Monographs. Oxford University Press, 2000.

\bibitem{karigiannis-2005-57}
S.~Karigiannis, {\it Deformations of {G}\_2 and {S}pin(7) {S}tructures on
  {M}anifolds},  {\em Canadian Journal of Mathematics} {\bf 57} (2005) 1012
  [\href{http://arXiv.org/abs/math/0301218}{{\tt math/0301218}}].

\bibitem{karigiannis-2007}
S.~Karigiannis, {\it Flows of ${G_2}$-{S}tructures, {I}},  {\em Q. J. Math.}
  {\bf 60} (2009), no.~4 487--522
  [\href{http://arXiv.org/abs/math/0702077}{{\tt math/0702077}}].

\bibitem{KarigiannisMcKayTsui}
S.~Karigiannis, B.~McKay and M.-P. Tsui, {\it {Soliton solutions for the
  Laplacian coflow of some $G_2$-structures with symmetry}},  {\em Differential
  Geom. Appl.} {\bf 30} (2012), no.~4 318--333
  [\href{http://arXiv.org/abs/1108.2192}{{\tt 1108.2192}}].

\bibitem{LotayWei1}
J.~D. Lotay and Y.~Wei, {\it Laplacian flow for closed {$G_2$} structures:
  {S}hi-type estimates, uniqueness and compactness},  {\em Geom. Funct. Anal.}
  {\bf 27} (2017), no.~1 165--233 [\href{http://arXiv.org/abs/1504.07367}{{\tt
  1504.07367}}].

\bibitem{LotayWei2}
J.~D. Lotay and Y.~Wei, {\it Laplacian flow for closed {$G_2$} structures: real
  analyticity},  {\em Communications in Analysis and Geometry} (2018)
  [\href{http://arXiv.org/abs/1601.04258}{{\tt 1601.04258}}]. in press.

\bibitem{LotayWei1a}
J.~D. Lotay and Y.~Wei, {\it Stability of torsion-free {$G_2$} structures along
  the {L}aplacian flow},  {\em Journal of Differential Geometry} (2018)
  [\href{http://arXiv.org/abs/1504.07771}{{\tt 1504.07771}}]. in press.

\bibitem{ManeroOtalVillacampa}
V.~Manero, A.~Otal and R.~Villacampa, {\it Solutions of the {L}aplacian flow
  and coflow of a {L}ocally {C}onformal {P}arallel {G2}-structure},
  \href{http://arXiv.org/abs/1711.08644}{{\tt 1711.08644}}.

\bibitem{MorganTianRicci}
J.~Morgan and G.~Tian, {\em Ricci flow and the {P}oincar\'e conjecture}, vol.~3
  of {\em Clay Mathematics Monographs}.
\newblock American Mathematical Society, Providence, RI, 2007.

\bibitem{PerelmanRicci}
G.~Perelman, {\it The entropy formula for the {R}icci flow and its geometric
  applications},  \href{http://arXiv.org/abs/math/0211159}{{\tt math/0211159}}.

\bibitem{ShiEstimate}
W.-X. Shi, {\it Deforming the metric on complete {R}iemannian manifolds},  {\em
  J. Differential Geom.} {\bf 30} (1989), no.~1 223--301.

\end{thebibliography}\endgroup

\end{document}